# Development of Testing Methodology in Mathematics Education in the Context of Digitalization

**Igor V. Orlovskyi[a] , Olena A. Tymoshenko[a,b]**

**[a]"Department of Mathematical Analysis and Probability Theory, NTUU Igor Sikorsky Kyiv Polytechnic Institute, Ukraine"**

**[b]"Department of Mathematics, University of Oslo, Norway"**

**Abstract**

Access to quality education remains a global challenge, particularly in crisis-affected regions. This study examines the decline in students' mathematical proficiency and proposes an innovative Moodle-based testing system that incorporates step-by-step solution verification and interactive exercises. Unlike traditional assessments, this approach ensures a more structured evaluation process, reducing student errors and improving feedback quality. Additionally, the study integrates Classical Test Theory to analyse test reliability, offering a novel perspective on the effectiveness of digital assessments. The proposed methodology improves assessment transparency and conceptual understanding, while Classical Test Theory is used to evaluate the reliability of multistage mathematical tasks.

## Introduction

In today's world, education should be accessible to everyone, regardless of location, economic status, or global crises. However, access to knowledge remains a privilege rather than a universal right. In developing countries, millions of children and adults face limitations due to poverty, lack of infrastructure, and social barriers. For instance, a study analysing the impact of conflict on education in Palestine found that conflict events significantly reduce high school students' probability of passing their final exams and being admitted to university (Brück, Di Maio, & Miaari, 2019).

Distance education has long been an established component of the educational landscape, supporting a wide range of learning contexts, including flexible access to education, lifelong learning, and geographically distributed instruction. It also plays a particularly important role in situations where traditional face-to-face education is disrupted.

The COVID-19 pandemic significantly expanded the scale of distance learning worldwide. School and university closures led to a rapid and often unplanned transition to online formats, revealing both the potential of distance education and a number of significant challenges. Similar changes were observed in undergraduate mathematics teaching, including linear algebra and calculus courses, where course design and delivery had to be adapted rapidly in response to the pandemic (see Alarfaj, O'Hagan, & Sangwin, 2022). This transition did not affect all educational contexts in the same way: research comparing online and traditional learning environments suggests that student progress may vary depending on the setting (Hurlbut, 2018). In particular, the shift posed notable difficulties in disciplines requiring practical experience (e.g., Akyildiz, 2021). More recently, other large-scale disruptions, such as armed conflicts, including the war in Ukraine and escalating tensions in Palestine, have further highlighted the importance of accessible and adaptable educational solutions. These situations have left many learners without stable access to education and have contributed to learning loss, as evidenced by reports on declining outcomes in language, reading, and mathematics (UNICEF, report, 2023).

Beyond immediate disruptions, young people who are forced to leave their home countries due to conflict must also adapt to new educational systems. This presents significant challenges, as seen in a study comparing the experiences of Ukrainian students learning mathematics in UK schools. The research highlights how differences in curricula, teaching methods, and assessment approaches create additional barriers to learning, further complicating students' academic progress. These challenges emphasize the need for adaptive and high-quality distance learning tools that can support displaced students in continuing their education effectively (see Proshkin & Foster, 2024).

In these circumstances, distance learning has become increasingly prominent. However, its effectiveness depends on the quality of tools, the accessibility of technologies, and methodologies adapted to different population groups.

This study presents an approach to integrating step-by-step solution verification in Moodle-based testing, ensuring a more structured and transparent assessment process. Unlike traditional assessment methods, this approach is intended to enhance interactivity and provide immediate feedback, which may help students identify and correct their mistakes during the assessment process. Furthermore, the research contributes to the field by analysing test reliability using Classical Test Theory (e.g., Crocker & Algina, 1986 or Wu, Tam, & Jen, 2016), offering insights into the consistency of this assessment approach.

This paper presents an approach for designing and evaluating stepwise mathematics assessment tasks implemented in Moodle, contributing to the study of distance mathematics education. Such tasks may enhance the accuracy and consistency of student knowledge assessment while also supporting conceptual understanding and spatial reasoning in mathematics learning environments (Schoenherr et al., 2024; Fahlgren et al., 2024).

The paper is organized as follows:

Section 1 discusses the role of mathematics in STEM and other disciplines, emphasizing the challenges students face in mastering mathematical concepts, particularly in distance learning settings. Section 2 outlines the problem context by presenting evidence of difficulties in students’ mathematical proficiency, including survey data and findings from international assessments, and considers the limitations of traditional assessment methods in such contexts. Section 3 presents the development of Moodle-based step-by-step testing and interactive exercises, demonstrating their potential to support more structured assessment and conceptual understanding. Section 4 analyzes test reliability based on Classical Test Theory, evaluating the consistency of the proposed assessment methods. Finally, the conclusion summarizes key findings and suggests directions for future research and implementation.

## 2. Statement of the Problem

Mathematics is an integral part of educational programs in STEM disciplines and requires continuous skill maintenance through regular practice. However, as practice shows, students often have significant knowledge gaps, which make it difficult to master the subject.

Survey data, which we collected in December 2024 in collaboration with the AdaMatS (Adaptive Mathematics Support in STEM Education) consortium, provide additional insight into this issue. The survey involved 120 university educators from 35 higher education institutions in Ukraine who teach first-year students in STEM programs, including both mathematics and professionally oriented courses. The study was designed to explore gaps in the mathematical preparation of first-year undergraduate students from the perspective of instructors. The questionnaire included closed, Likert-scale, semi-closed, and open-ended questions, allowing respondents to evaluate students' knowledge levels, identify problematic areas, and reflect on possible causes of these gaps. In this context, "successful mastery" refers to instructors' perceptions of students' ability to confidently apply core mathematical concepts and methods in solving standard and applied problems. "Sufficient understanding" reflects a level at which students demonstrate partial comprehension but may experience difficulties in applying knowledge independently, while "adequate academic performance" refers to achieving minimum course requirements, typically reflected in passing grades. The questionnaire used in this study is publicly available (AdaMatS Consortium, 2024).

The survey results indicate the following patterns in lecturers' perceptions:

- only 8.4% of lecturers believe that 80% of students successfully master mathematical disciplines;
- 33.3% of respondents estimate that only 40% of students achieve sufficient understanding, while 29.2% believe that only 20% of students demonstrate adequate academic performance;
- more than a third of lecturers (34.2%) report that students' mathematical knowledge is declining, and nearly half (48.3%) believe that a decrease is more likely than an improvement.

These findings reflect instructors' perceptions of students' mathematical preparation rather than direct measurements of learning outcomes. To place these perceptions in a broader context, it is useful to consider evidence from international studies reporting declines in mathematics assessment outcomes.

For instance, the paper *"Should Mathematicians Worry with PISA and TIMSS Math Results?"* (Nuno & Maroco, 2024) argues that declines in mathematical literacy among European students may pose serious challenges for scientific progress, economic development, and education quality. Similarly, the article *"Mathematics Test Scores in Some Countries Have Been Dropping for Years, Even as the Subject Grows in Importance*" (Samuels, 2024) reports declining scores on international mathematics assessments in countries such as France, Germany, and Finland. However, such declines should be interpreted cautiously, as they may also reflect differences in curricula, test content, or educational expectations. The study *"Changes in Achievement on PISA: The Case of Ireland and the United Kingdom"* (Cosgrove & Cartwright, 2014) likewise identifies downward trends in PISA results in these countries.

Taken together, these studies provide external evidence of challenges in mathematics achievement, while our survey captures instructors' perceptions of these challenges in the context of first-year STEM education.

### *2.1 Challenges in Mathematics Education and Assessment*

In mathematics education, a high level of proficiency involves the development of both conceptual understanding and problem-solving skills, which are often supported through interaction with instructors and structured learning activities. Previous research highlights the importance of interactivity in learning, emphasizing its role in supporting students' cognitive engagement and understanding in digital environments (Moreno & Mayer, 2007; Cairncross & Mannion, 2001). In situations such as wars, pandemics, and other disruptions, in-person instruction may become unavailable, necessitating a transition to remote learning formats. Under such conditions, students may experience additional difficulties related to limited access to instructional support, as well as technical constraints such as low internet speed, power outages, and restricted access to digital resources (Deiser & Reiss, 2014; Lindenbauer et al., 2024).

While asynchronous online learning environments provide flexible access to educational content, they may present challenges for students who require timely feedback and guidance when learning complex mathematical concepts. In addition, crisis conditions may increase academic and mathematics-related anxiety, which can reduce students' concentration and motivation and further hinder mathematics learning (Caviola et al., 2021).

In this context, traditional teaching and assessment approaches may require adaptation to better support students' learning in distance and disrupted learning environments.

### *2.2 The Need for Modernized Testing Approaches*

In mathematics education, these limitations become even more pronounced in digital assessment environments. Many traditional test formats, such as multiple-choice questions, may allow students to arrive at the correct answer through elimination or guessing, which can limit the extent to which such assessments capture students' reasoning processes and deeper conceptual understanding. Moreover, the rise of artificial intelligence tools like ChatGPT and Wolfram Alpha enables students to access ready-made solutions almost instantly, diminishing the effectiveness of standard testing methods. Recent research on e-assessment in undergraduate mathematics also emphasizes the need to align assessment design with the specific demands of mathematical learning. In particular, a collaboratively developed research agenda identifies key challenges related to feedback, student interaction, task design, the affordances of e-assessment tools, and the assessment of mathematical skills in digital environments (Kinnear et al., 2024). These developments underscore the need for assessment strategies that not only verify the correctness of answers but also evaluate the reasoning behind them.

Research on assessment and feedback suggests that testing can fulfil both learning-oriented and evaluative functions (Shute, 2008; Sangwin, 2013):

- *Educational* - helping students develop problem-solving skills rather than merely testing their ability to recall facts.

- *Evaluative* - objectively measuring students' knowledge and understanding in a consistent and transparent way.

To achieve these objectives, modern testing methods may incorporate the following principles (Sangwin, 2013; Gallaun et al., 2022; Kinnear et al., 2024; Fahlgren et al., 2024):

- assess not only the final result but also the step-by-step reasoning process;
- include structured problem-solving tasks, in addition to selected-response formats such as multiple-choice questions;
- adapt to students' proficiency levels by offering varying degrees of difficulty;
- implement strategies that reduce the risk of over-reliance on guessing, answer sharing, or AI-generated solutions.

In particular, recent studies in mathematics e-assessment emphasize high-quality task design, follow-up error handling, randomization, and meaningful feedback as key conditions for effective digitally supported mathematical problem solving (Gallaun et al., 2022; Kinnear et al., 2024). These principles directly inform the Moodle-based stepwise task design presented in the next section.

Another challenge in the modernization of assessment methods is the limited time and expertise available to educators for analysing test quality. Designing effective tests may require digital literacy, pedagogical expertise, and access to appropriate technical support, as well as a deep understanding of assessment reliability and validity (Husain, 2021; Pastore, 2023). International organizations, such as the International Test Commission, provide standards intended to ensure fair and appropriate assessment practices (International Test Commission, 2006). However, successfully implementing these guidelines necessitates additional teacher training and institutional support.

In the context examined in this study, improving mathematics assessment may require changes in both teaching methodologies and assessment techniques. One promising approach is the integration of adaptive learning strategies, such as didactic transposition (El Fadel, 2024), understood as the adaptation of complex disciplinary knowledge into forms suitable for teaching and learning while preserving academic rigor. Additionally, Shulman's framework for pedagogical content knowledge (Shulman, 1986), which highlights the connection between subject knowledge and effective ways of teaching and assessing it, can be applied to strengthen students' understanding of mathematical structures and problem-solving techniques.

A key aspect of this modernization is the development of more diverse and flexible test formats, including selected-response items, coded-response tasks, step-by-step problem-solving exercises, and interactive tasks implemented within digital platforms. In addition to traditional answer-selection tasks, structured problem-solving exercises can support a more comprehensive evaluation of students' conceptual understanding. Adaptive testing, which adjusts question difficulty based on responses, may also support a more precise diagnosis of students' knowledge gaps, while requiring coded responses and explanations of problem-solving steps may reduce the likelihood of guessing or reliance on external resources such as AI tools or online solvers.

For such assessments to be both effective and scalable, digital platforms should support features aligned with the specific requirements of mathematics education. These include the ability to diagnose gaps in students’ knowledge and understanding, provide timely and detailed feedback, and support step-by-step problem solving. Research on mathematics teaching, learning, and assessment in the digital age highlights the growing importance of digital tools in supporting these processes (Weigand et al., 2024). Interactive elements such as visualizations and simulations can assist in explaining abstract concepts (Schoenherr et al., 2024), while adaptive task design can enable a gradual increase in task complexity. In addition, automatically generated task variants and question banks can support more consistent assessment conditions and reduce the likelihood of answer sharing. The effective implementation of such digital assessment practices also depends on teachers' assessment-related competencies in digital environments (Estaji et al., 2024). More broadly, digital and blended learning environments, including institutional learning management systems and course-specific digital tools, can provide students with supplementary materials, preparatory modules, and practice tasks that help bridge gaps in mathematical knowledge and support the transition to more advanced mathematical study.

Taken together, these considerations highlight the need for more structured and interactive assessment formats in mathematics education, particularly in digital learning environments. In this context, step-by-step testing supported by digital tools may provide a more detailed picture of students' reasoning and problem-solving processes. The following section examines the practical implementation of this approach within the Moodle platform.

Within the broader field of mathematics e-assessment, specialised systems such as STACK and Numbas are widely used because they support symbolic validation, randomized task generation, and adaptive feedback (Sangwin, 2013; Kinnear et al., 2024; Fahlgren et al., 2024). In contrast, the present study explores how similar pedagogical principles can be implemented directly within a standard Moodle environment using its native tools, thereby lowering the entry threshold for both instructors and students already familiar with the LMS.

## 3. Enhancing Assessment Accuracy with Step-by-Step Testing in Moodle

Building on the challenges outlined in the previous sections, this section presents a practical approach to enhancing mathematics assessment through step-by-step testing implemented in Moodle. At Igor Sikorsky Kyiv Polytechnic Institute, the Sikorsky distance learning platform has been developed based on Moodle (Modular Object-Oriented Dynamic Learning Environment). Originally developed in Australia by WebCT employees, Moodle is an open-source learning management system (LMS) that provides extensive tools for designing and managing online assessments. Previous studies have also highlighted its potential for adaptive learning and the implementation of adaptive mechanisms in higher education settings (Karagiannis & Satratzemi, 2020; Morze et al., 2021). The Sikorsky Moodle-based platform allows educators to create and edit course content with integrated text, graphics, animations, and multimedia elements. It also facilitates the customization of test tasks, provides detailed course statistics for monitoring student progress, and supports discussions through an integrated forum. Additionally, instructors can configure test settings, including the number of attempts, time limits, question and answer randomization, and flexible grading based on question complexity.

The Sikorsky platform hosts numerous Moodle-based online courses, providing students with access to educational materials and diverse assessments to aid their learning and self-evaluation in mathematics. However, further improvements are required to develop interactive tools that extend beyond traditional test formats, ensuring a more engaging and effective digital learning experience.

Moodle offers various tools for online assessment, such as multiple-choice tests, matching exercises, and drag-and-drop tasks. However, these standard formats often fail to fully assess students' understanding of mathematical concepts, as they focus primarily on selecting the correct answer rather than demonstrating the underlying thought process.

To address these limitations, we have developed specialized tests using Moodle's built-in tools. These tests enhance assessment validity by evaluating students' logical reasoning rather than just their final answers. Constructed using HTML and LaTeX, they ensure clear and precise mathematical formatting.

### *3.1 Determinant-Based Equation Solving*

This task is later used in Section 4 as a representative example of a compact stepwise mathematical task (later denoted as Task C) to demonstrate the full process of psychometric analysis within the proposed methodology.

Consider the following problem, where students must solve the equation:

$$\begin{vmatrix} x+4 & 5 & 3 \\ -5 & x-6 & -2 \\ 1 & 1 & 1 \end{vmatrix} = 0.$$

The test is structured to verify both the determinant expansion and the final solution:

**Step 1: Compute the Determinant**

Students expand the determinant to derive a quadratic equation. The test ensures that the coefficients align with the expected form

$$\text{Determinant} = (\ 1)x^2 + (-3)x + (2) = 0.$$

Validation: Moodle checks the individual coefficients $(1, -3, 2)$ to ensure accurate intermediate calculations.

**Step 2: Solve the Quadratic Equation**

Students then solve

$$x^2 - 3x + 2 = 0$$

and provide the roots in ascending order: $x_1 = 1, x_2 = 2$.

Validation: The system verifies both roots, ensuring logical correctness.

The test is implemented using embedded LaTeX for mathematical expressions and numerical input validation for stepwise grading. Figure 1 illustrates the structure of a Moodle-based question designed to assess determinant evaluation and quadratic equation solving.

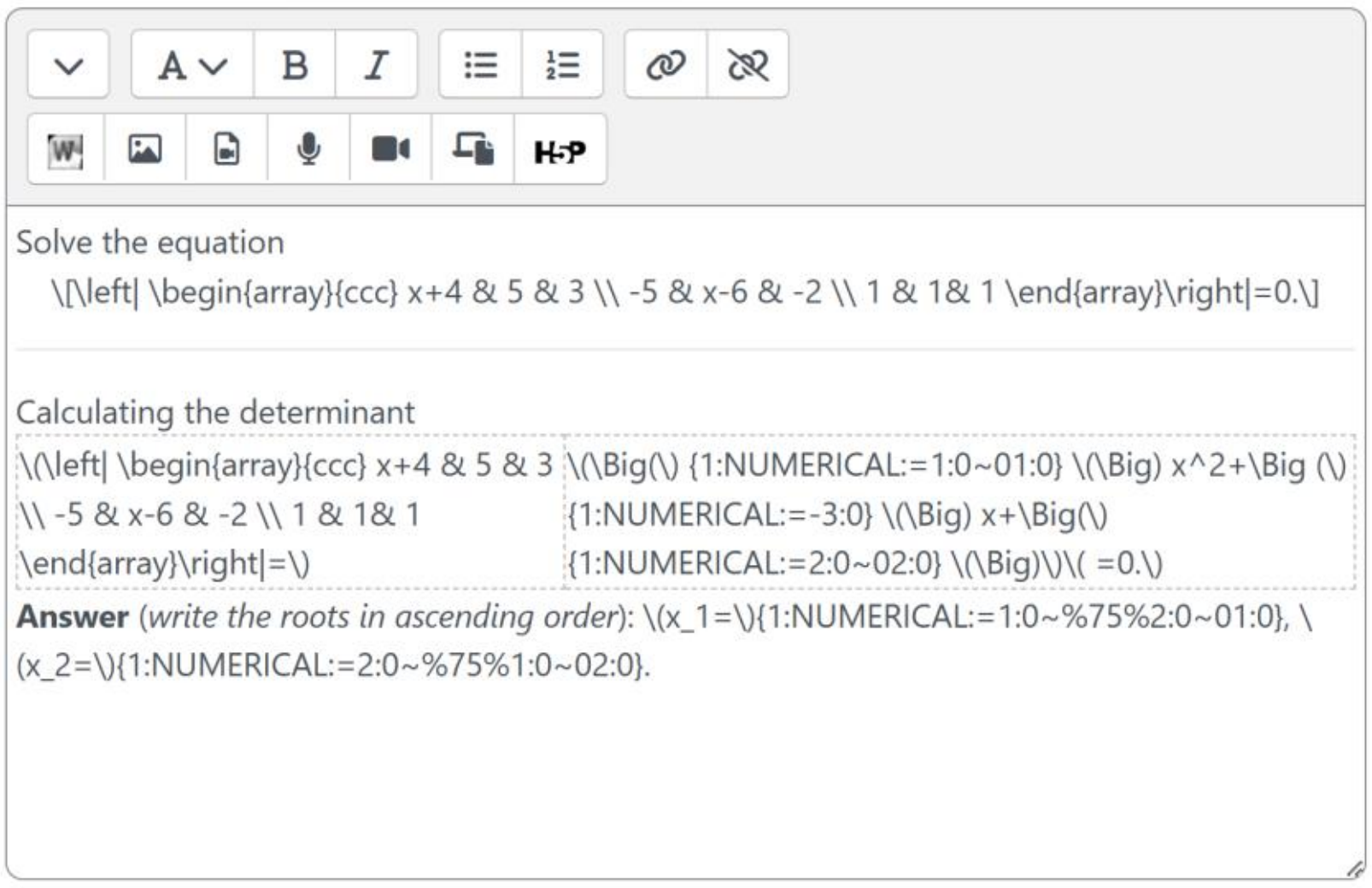

**Figure 1:** Moodle question implementation using LaTeX, displaying determinant evaluation, quadratic equation derivation, and root finding for $x^2 - 3x + 2 = 0$

Figure 2 presents the interactive Moodle test interface for solving matrix determinant problems. The structured layout helps students navigate the logical workflow, enabling automatic assessment of both intermediate calculations and final solutions. The interface uses visual aids, such as bullet points for coefficients and boxed answer fields, to enhance clarity and guide learners through the problem-solving process.

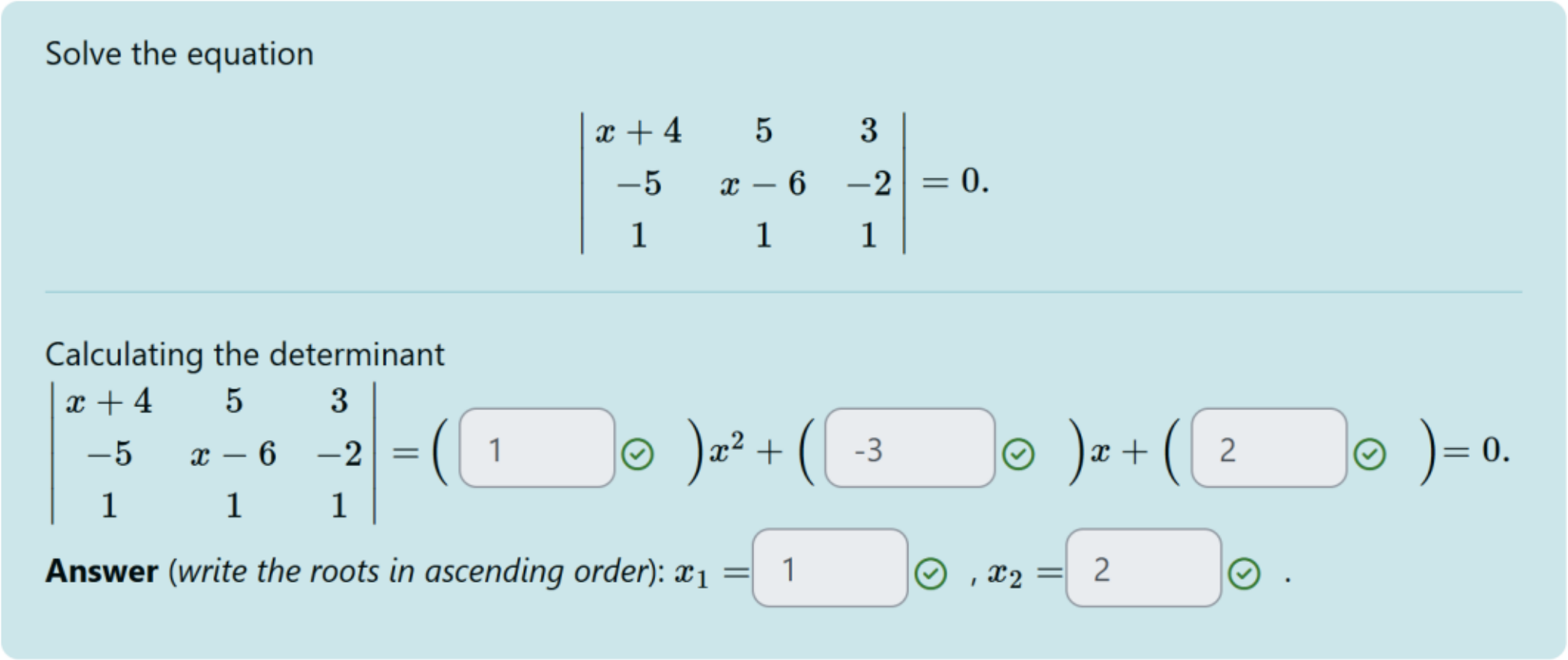

**Figure 2:** Moodle test interface for solving , showing stepwise evaluation and answer validation.

This testing methodology provides several key benefits. By requiring students to document their reasoning, it allows instructors to identify specific errors in the problem-solving process, such as mistakes in determinant expansion or algebraic transformations, and to provide targeted feedback. Additionally, partial credit can be awarded when intermediate steps are correct, ensuring a fairer and more comprehensive assessment focused on students' logical processes rather than just their final answers. Due to its clear internal sequence and moderate structural complexity, this task also serves as a convenient representative example for the detailed reliability analysis presented in Section 4.

*3.2 Interactive Integration Exercise*

The benefits of this structured assessment approach extend to other interactive learning tasks, such as the integration-based drag-and-drop exercise illustrated in Figure 3.

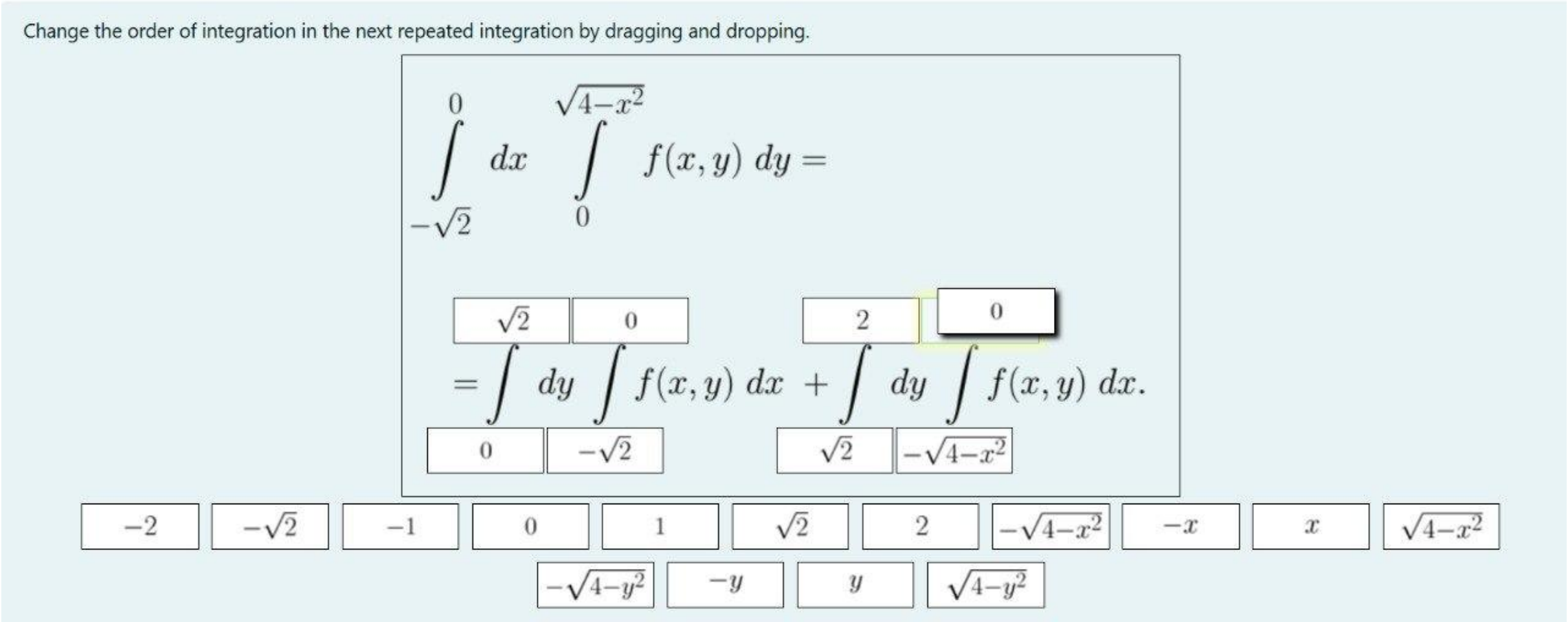

**Figure 3:** Interactive drag-and-drop exercise for determining integration bounds, reinforcing conceptual understanding

In this activity, students must correctly place integration limits by dragging the appropriate answer blocks into the corresponding empty slots. This process allows them to construct the integration bounds step by step, visually reinforcing the structure of the integral. By assembling the limits like a puzzle, students develop a deeper understanding of the order of integration and its geometric meaning.

The system provides instant feedback, helping to identify and correct common mistakes. For example, students may misplace negative and positive limits or incorrectly assign values to the wrong variables. If an error occurs, the system highlights the mistake and offers guidance. Even if the final answer is not entirely correct, partial credit is awarded for properly positioned elements.

Another common approach to assessing students' understanding of integration order is through multiple-choice questions, as shown in Figure 4.

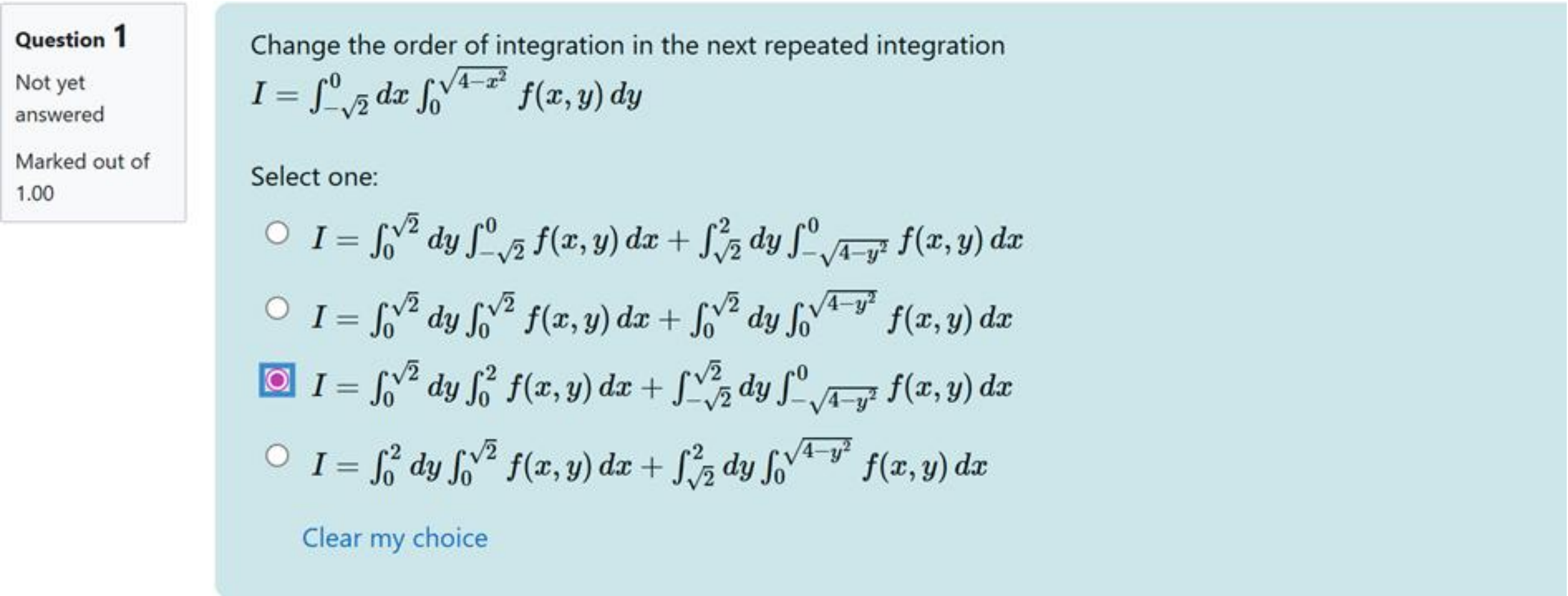

**Figure 4:** Multiple-choice exercise for determining integration bounds

The multiple-choice format allows for quick automatic grading, making it suitable for large-scale assessments. It provides clear and unambiguous results since students either select the correct answer or they don't. Additionally, it is easy to implement and does not require complex technical resources. However, it has some limitations, such as a lower degree of active engagement, the possibility of guessing the correct answer without understanding the concept, and the risk of students focusing on recall rather than analytical reasoning.

On the other hand, the drag-and-drop format offers a more interactive and engaging experience. By requiring students to actively construct the integral step by step, it promotes a deeper understanding of the topic. This approach may support structural understanding of the integral, as students must correctly organize and place its components. Unlike multiple-choice tests, where a correct answer can be chosen randomly, the drag-and-drop method reduces the chances of guessing and allows students to see their mistakes immediately. However, this format also presents several challenges, including the need for technical implementation, longer completion time, potential accessibility barriers for students with motor control difficulties, and the need for well-designed feedback to ensure that students understand their errors.

Overall, if the goal is to quickly assess students' knowledge of integration order, the multiple-choice format is more efficient. However, if the objective is to promote conceptual understanding and careful reasoning about integration limits, the drag-and-drop format may offer additional pedagogical value. The best approach would be to combine both methods: using interactive exercises to build understanding and reinforcing the knowledge with multiple-choice tests.

This interactive format may help students engage more actively with the structure of the integral. When combined with appropriate visual representations and immediate feedback, it may support both technical proficiency and conceptual understanding.

### *3.3 Interactive Drag-and-Drop Task for Matching Correct Answers to an Image*

Another example of interactive testing is the drag-and-drop task shown in Figure 5. Originally, this type of test in Moodle was designed for humanities disciplines, but we adapted it for mathematics, which is a significant step forward from a pedagogical perspective.

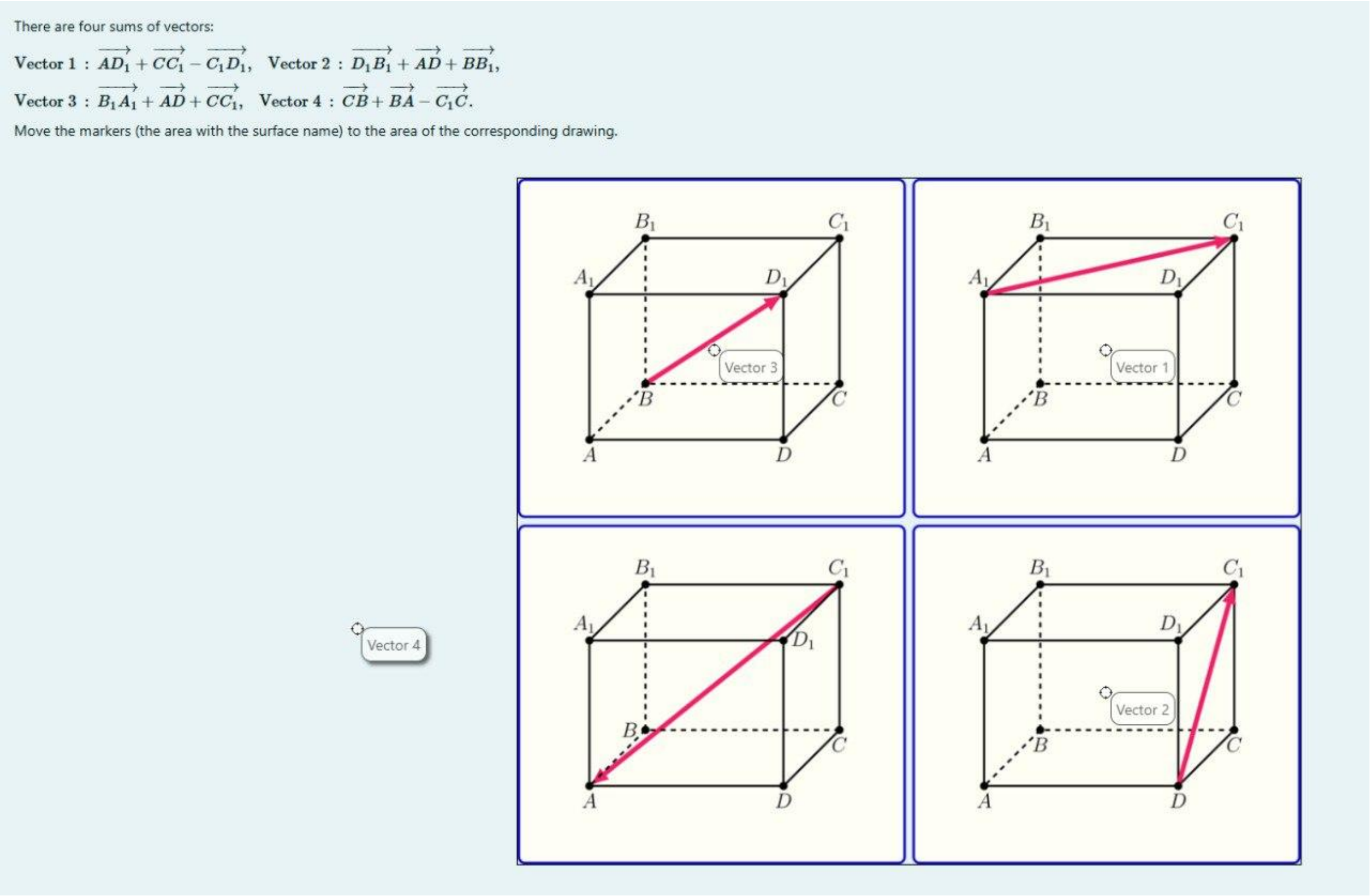

**Figure 5:** Interactive drag-and-drop exercise for matching vector sums with their corresponding graphical representations, reinforcing spatial reasoning and conceptual understanding.

Such drag-and-drop visual tasks are widely reported to support visual-spatial reasoning and increase student engagement in mathematics learning contexts (Schoenherr et al., 2024; Fahlgren et al., 2024).

In this exercise, students match vector transformations with their corresponding graphical representations by dragging elements to the appropriate positions. This interactive method reinforces spatial reasoning, enhances conceptual understanding, and bridges the gap between algebraic expressions and geometric interpretations.

Research in cognitive psychology supports the effectiveness of this approach, demonstrating that hands-on manipulation of digital objects improves memory retention and comprehension of abstract concepts. Additionally, transforming routine problem-solving into an engaging activity increases student motivation. Automated grading and instant feedback make this tool particularly valuable for large-scale education, streamlining the assessment process while ensuring accurate evaluation.

Regular integration of such exercises in vector algebra courses significantly enhances students' accuracy in vector operations, accelerates problem-solving, and strengthens their ability to visualize mathematical objects in space. This adaptability extends to various difficulty levels, from fundamental 2D vector operations to complex 3D transformations, making the tool suitable for both introductory and advanced mathematics.

### *3.4. STACK-Based Symbolic Assessment with Parsons problems blocks*

Building on the stepwise and interactive formats presented above, STACK-based tasks can further extend symbolic mathematics assessment within Moodle. Unlike standard drag-and-drop or numerical-response tasks, STACK enables automatic verification of algebraic equivalence and supports mathematically valid alternative answer forms.

In the present implementation, the **Parsons library** was used to construct interactive Parsons-style tasks in which students drag and arrange symbolic expressions into the correct positions. An illustrative example is a calculus task in which students are required to match each function $y(x)$, defined for $x > 0$, with its first and second derivatives. Figure 6.(*a*) presents a STACK-based exercise involving power, exponential, and logarithmic functions.

To increase variability and reduce answer memorization, the exponent in the power function, the coefficient multiplying $x$ in the exponential function, the coefficient inside the logarithm argument, and the logarithm base are generated randomly for each task instance. This parameterized generation allows the same assessment template to produce a large number of mathematically equivalent but numerically distinct task variants.

The task was implemented using the Parsons interaction model, where students construct the full derivative structure step by step by placing symbolic expressions into the rows corresponding to $y$, $y'$, and $y''$. Such a format supports procedural fluency in differentiation while reinforcing the structural relationships between a function and its higher-order derivatives.

An additional pedagogical advantage of STACK is the automatic generation of detailed answer feedback. After submission, the system produces a formatted solution table displaying the function together with its correct first and second derivatives (Figure 6.(*b*)). This immediate symbolic feedback supports error analysis, self-correction, and deeper understanding of derivative patterns.

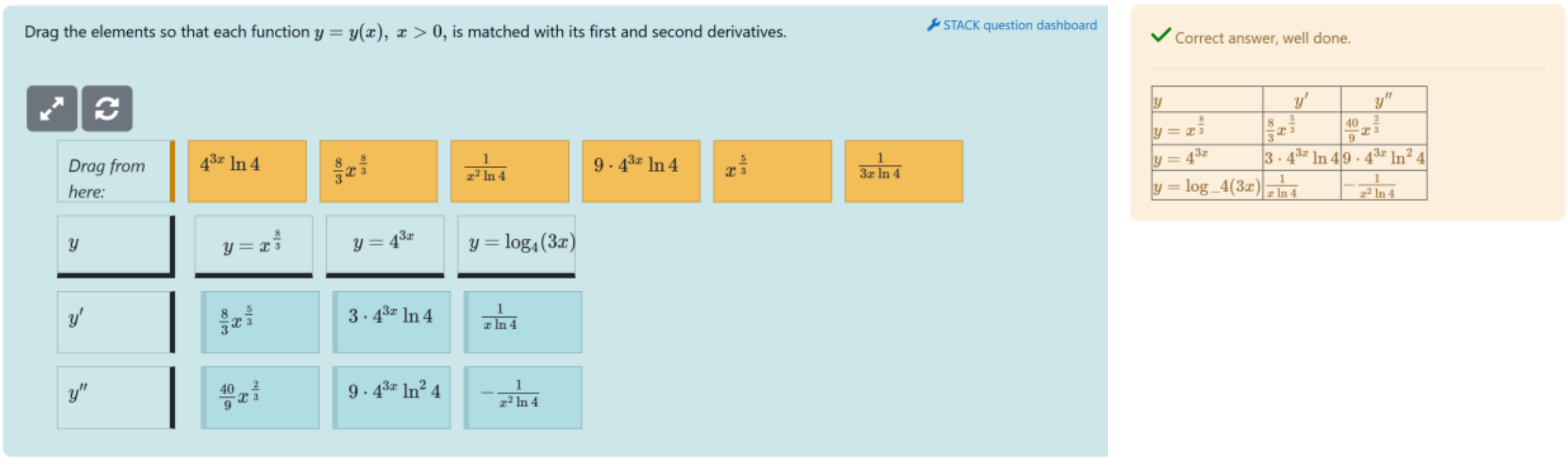

(*a*) (*b*)

**Figure 6:** STACK-based Parsons task for matching a function with its first and second derivatives, with randomly generated parameters and automatically generated STACK feedback table displaying the correct symbolic derivatives after task completion.

Such tasks naturally extend the stepwise assessment principle developed in this study. In addition to immediate feedback, STACK provides symbolic validation and randomized parameter generation, making it particularly suitable for calculus, algebraic transformations, and multistage mathematical reasoning.

In the future, we plan to enhance this approach by introducing more complex mathematical problems and adaptive tests that automatically adjust to each student's proficiency level. A particularly promising direction is the deeper integration of STACK into stepwise assessment formats, enabling symbolic validation not only of final answers but also of intermediate transformations and logically connected micro-steps. This will allow more flexible implementation of multistage algebraic and calculus tasks, including randomized parameter generation. Such extensions may improve student engagement and make multistage symbolic tasks more effective in mathematics learning, as suggested by recent studies on computer algebra-based assessment environments (Sangwin, 2013; Kinnear et al., 2024).

## 4. Analysis of Test Reliability Based on Classical Test Theory

The development of stepwise tests is a challenging task, as it requires an individual approach to each problem and structuring the solution into a sequence of logical steps that can be automatically evaluated within a digital platform. However, the implementation of new assessment formats requires not only their technical realization but also a rigorous psychometric evaluation of their quality and reliability.

To assess the effectiveness and internal consistency of the proposed testing methodology, an empirical analysis was conducted using Classical Test Theory (Crocker & Algina, 1986). Before full implementation, the developed tasks were piloted with a group of students, and the resulting data were analysed to evaluate the reliability of both the entire four-task test battery and its individual components.

As an example of such an analysis, we focus on a task related to finding the determinant of a matrix (see Section 3.1) and examine the reliability of the entire test, which includes this task.

The test was designed to assess the knowledge of engineering students during the first semester of a three-semester course in higher mathematics. A total of 117 students from five study groups of the Faculty of Informatics and Computer Engineering at Igor Sikorsky Kyiv Polytechnic Institute participated in the study during the 2025-2026 academic year. The course carries 17 ECTS credits.

The test consisted of four tasks:

Task A: finding a product of two matrices (see Figure 7.(*a*));

Task B: finding a matrix polynomial (see Figure 7.(*b*))

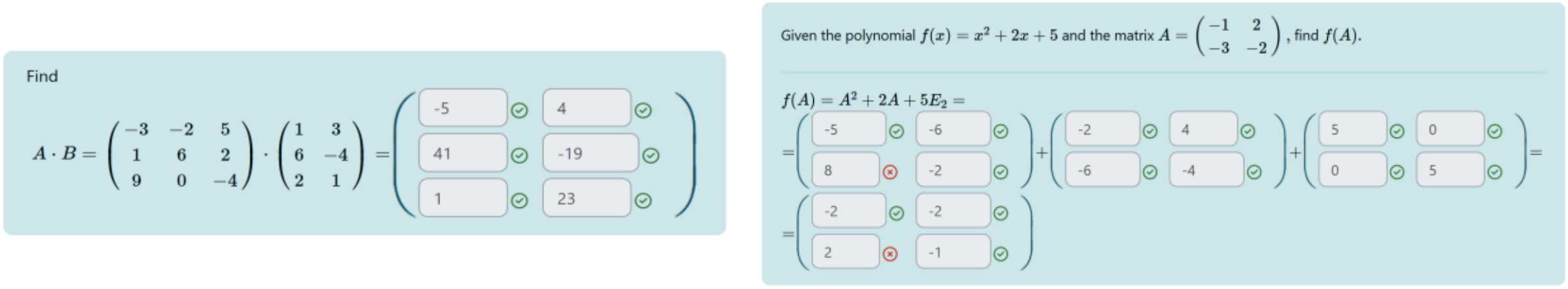

(*a*) (*b*)

**Figure 7.**

Task C: determining the values of $x$ that make the matrix determinant equal to zero. The solution involves forming and solving a quadratic equation (see Figure 2);

Task D: Calculating the determinant of a fourth-order matrix by reducing it to triangular form (see Figure 8).

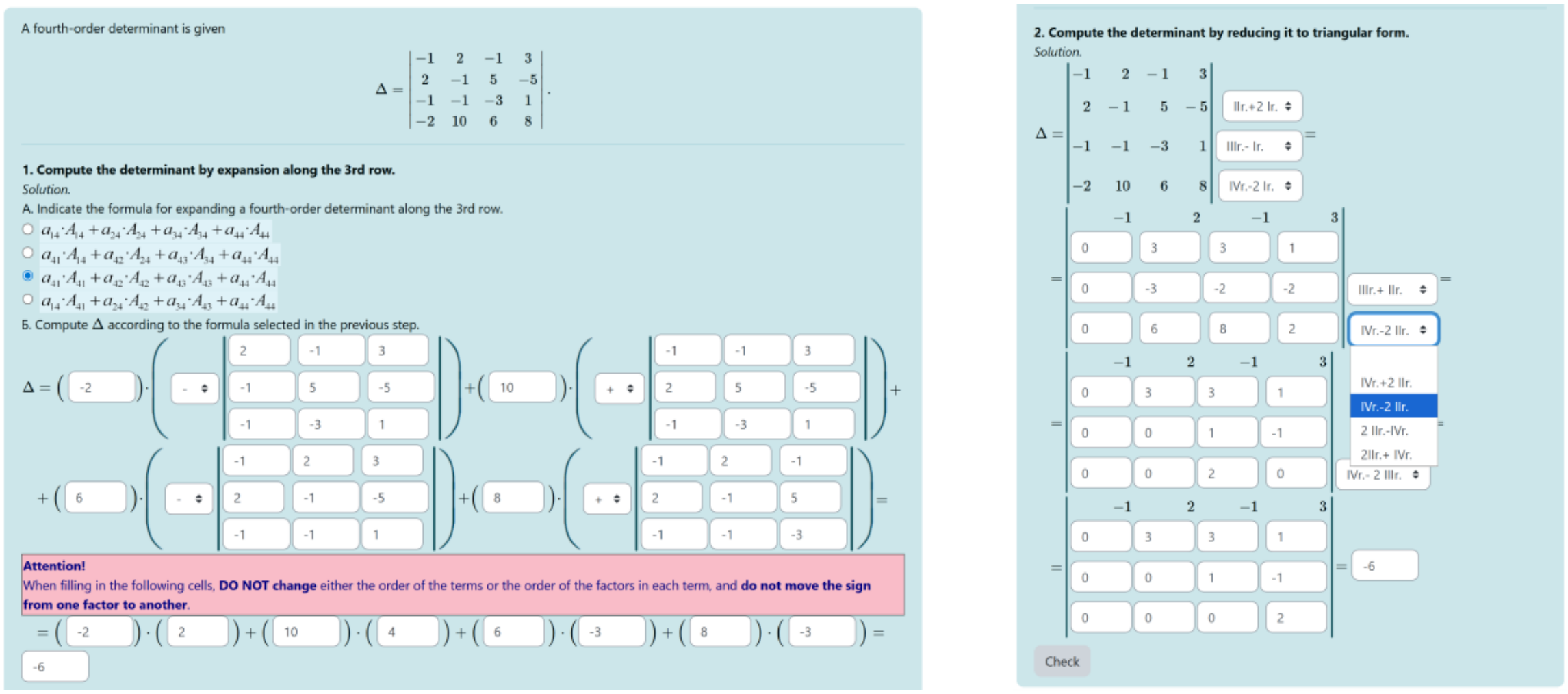

**Figure 8.**

Building on the representative determinant-based task introduced in Section 3.1, we next examine the reliability of the full four-task battery and then use Task C as an illustrative example of detailed task-level psychometric analysis. Its structure is sufficiently rich to illustrate the full process of task decomposition, Moodle-based implementation, and subsequent psychometric evaluation, while remaining compact enough for transparent step-by-step interpretation. The task consists of two main stages: calculating the determinant to obtain the coefficients of a quadratic equation and then determining its roots.

The following subsection describes the statistical methodology used to evaluate score distributions, scale reliability, and item-level characteristics of the proposed stepwise assessment format.

### *4.1 Research Methodology*

Within the framework of Classical Test Theory, the observed test scores were analysed using descriptive statistics, reliability measures, and item-level indicators.

The analysis begins with the compilation of a dataset of students' test scores. The anonymized normalized score dataset used for the Classical Test Theory analysis, together with task-level scoring tables and supplementary figures, is publicly available in the Zenodo repository (DOI: 10.5281/zenodo.19480465). For each student, individual task scores $x_i$ were recorded, where $i = 1, \dots, N$ and $N$ denotes the number of students, and all task scores were normalized to the interval $[0; 1]$ to ensure comparability across tasks. The overall score was computed as a weighted sum of the task scores using weights shown in Table 1, and was subsequently normalized to the same interval $[0; 1]$.

To analyse the distribution of test results, several key descriptive statistical measures were calculated. The sample mean $\bar{x}$ represents the average score

$$\bar{x} = \frac{1}{N}\sum_{i=1}^{N} x_i \,.$$

In addition, the median and mode were used as robust indicators of central tendency.

The variability of individual scores was measured by the standard deviation

$$S_x = \sqrt{\frac{1}{N-1}\sum_{i=1}^{N}(x_i - \bar{x})^2}\,.$$

The shape of the score distribution was analysed using the skewness coefficient

$$Sk = \frac{\frac{1}{N}\sum_{i=1}^{N}(x_i - \bar{x})^3}{S_x^3}$$

which indicates the degree of asymmetry of the distribution, and the kurtosis coefficient

$$K = \frac{\frac{1}{N}\sum_{i=1}^{N}(x_i - \bar{x})^4}{S_x^4} - 3$$

which characterizes the concentration of values around the mean.

In addition to descriptive statistics, several indicators of test reliability and item quality were evaluated within the framework of Classical Test Theory.

To evaluate the internal consistency and psychometric properties of the test, several reliability indicators were calculated, including Cronbach's alpha and McDonald's omega.

Cronbach's $\alpha$ measures internal consistency, indicating how well the test questions measure the same concept. It is calculated using the formula:

$$\alpha = \frac{K}{K-1}\left(1 - \frac{\Sigma\sigma_i^2}{\sigma_t^2}\right),$$

where *K* is the number of test items, $\sigma_i^2$ is the variance of each question, and $\sigma_t^2$ is the variance of the total test score.

McDonald's $\omega$ is considered a more precise measure, as it accounts for differences in item variance. It is given by:

$$\omega = \frac{(\Sigma\lambda_i)^2}{\Sigma\lambda_i^2 + \Sigma\sigma_i^2},$$

where $\lambda_i$ represents factor loadings (weight of each question in the model) and $\sigma_i^2$ is the specific variance of each question.

Finally, item-level characteristics were analysed to assess the contribution of individual tasks to the overall test performance. Besides analysing correlations between tasks, the following item-level characteristics were considered.

At the item level, the relationship between the score of an individual task *X* and the total score of the remaining test tasks *R* was evaluated using the corrected item-rest correlation coefficient

$$r_{XR} = \frac{\sum_{i=1}^{N}(x_i - \bar{x})(R_i - \bar{R})}{\sqrt{\sum_{i=1}^{N}(x_i - \bar{x})^2 \cdot \sum_{i=1}^{N}(R_i - \bar{R})^2}},$$

where $R_i$ is the total score of the remaining tasks for student $i$ (excluding the analysed task), and $\bar{R}$ is the corresponding mean rest score.

Additionally, the difficulty of an individual task was estimated by the item difficulty index

$$p = \frac{\bar{x}}{x_{max}},$$

where $x_{max}$ is the maximum possible score.

*4.2 Results*

The descriptive statistics for the four stepwise tasks (Test A–D) and the normalized weighted overall score are presented in Table 1.

**Table 1.** ***Descriptive Statistics***

| | Test A | Test B | Test C | Test D | Score |
|---|---|---|---|---|---|
| **Weight** | 1 | 1 | 1 | 3 | - |
| **Mean** | 0.974 | 0.947 | 0.930 | 0.857 | 0.904 |
| **Median** | 1.00 | 1.00 | 1.00 | 0.940 | 0.940 |
| **Mode** | 1.00 | 1.00 | 1.00 | 1.00 | 1.00 |
| **Standard deviation** | 0.0985 | 0.140 | 0.172 | 0.202 | 0.117 |
| **Skewness** | -4.17 | -3.23 | -3.68 | -2.19 | -1.94 |
| **Kurtosis** | 17.4 | 10.5 | 15.4 | 5.04 | 4.29 |
| **Item-rest correlation** | 0.3022 | 0.0860 | 0.0782 | 0.1683 | - |
| **Shapiro-Wilk W** | 0.285 | 0.435 | 0.462 | 0.716 | 0.780 |
| **Shapiro-Wilk p** | <.001 | <.001 | <.001 | <.001 | <.001 |

The obtained results demonstrate a pronounced concentration of scores near the maximum values across all four tests. The mean scores for Test A, Test B, Test C, and Test D were 0.974, 0.947, 0.930, and 0.857, respectively. The normalized weighted overall score had a mean value of 0.904, indicating that most students successfully completed the majority of the stepwise tasks. This pattern is also clearly visible in the histogram of the normalized overall score (Figure 9.$(a)$), where the majority of observations are concentrated near the upper boundary of the scale.

The comparison of mean, median, and mode values reveals a clear ceiling effect, with the median and mode reaching or approaching the maximum possible value for all tests. This interpretation is reinforced by the strongly negative skewness coefficients, ranging from $-4.17$ to $-1.94$, which indicate that the score distributions are heavily concentrated in the upper range of the scale. The high kurtosis values (from 4.29 to 17.4) further confirm the peaked character of the distributions and the substantial proportion of near-perfect responses.

**Table 2.** ***Scale Reliability Statistics***

| | **Mean** | **SD** | **Cronbach's α** | **McDonald's ω** |
|---|---|---|---|---|
| **scale** | 0.927 | 0.0886 | 0.275 | 0.455 |

The reliability analysis revealed low internal consistency of the four-task test battery. The obtained Cronbach’s alpha was 0.275, while McDonald’s omega reached 0.455, both of which remain below commonly accepted thresholds for educational assessments. These values indicate that the four tasks are better understood as related but distinct components of mathematical competence rather than as a single unified scale.

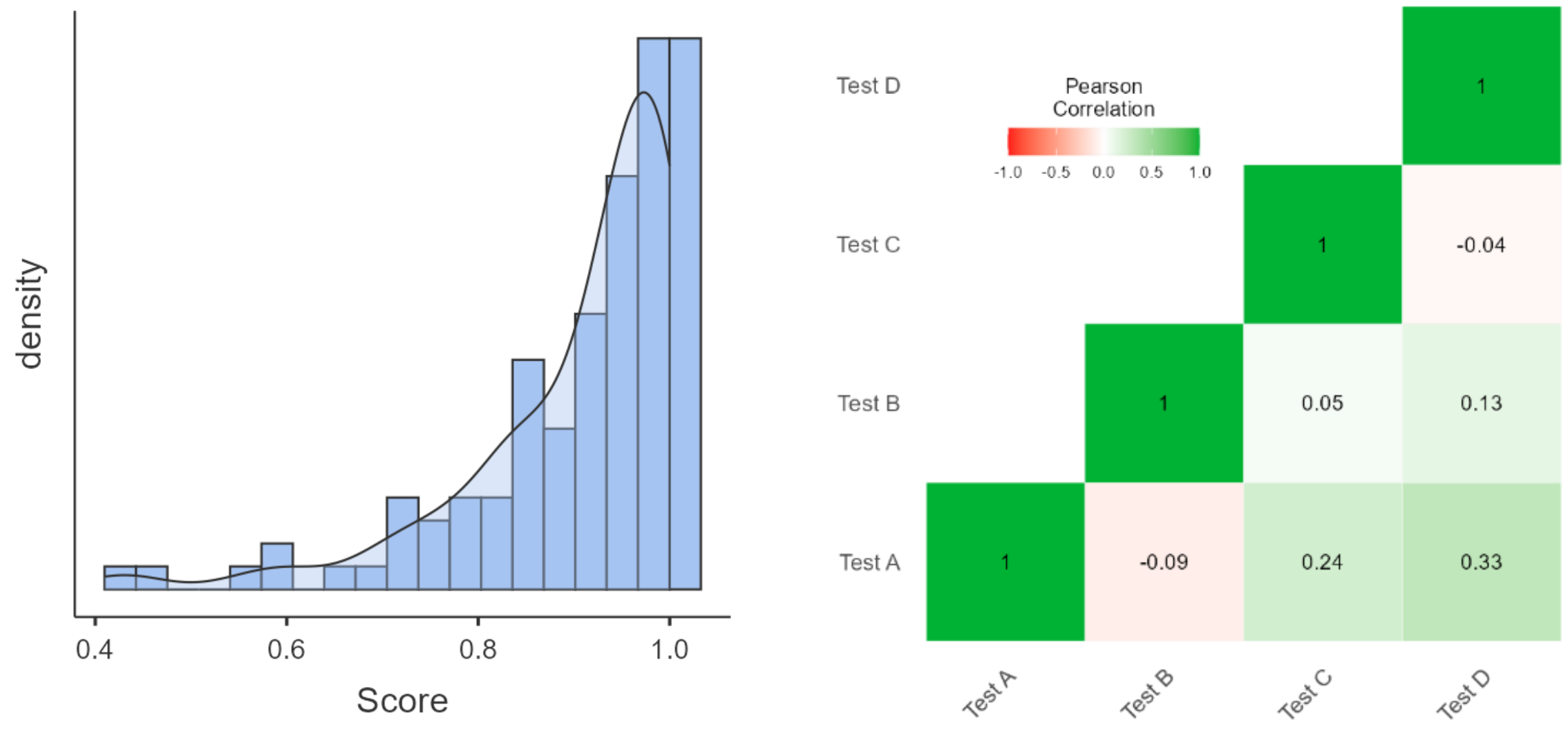

($a$) Histogram (b) Correlation heatmap

**Figure 9:** Histogram and correlation heatmap between task scores.

This interpretation is also visible in the correlation heatmap (Figure 9.($b$)). Inter-test Pearson correlations were generally weak, ranging from $-0.09$ to $0.33$. In particular, the slightly negative relationship between Test C and Test D ($r = -0.04$) and the near-zero correlation between Test B and Test C ($r = 0.05$) indicate that different tests may require partially distinct mathematical operations and reasoning strategies. Since the four tests correspond to different stages of matrix- and determinant-related problem solving, the observed weak relationships may reflect the heterogeneous cognitive structure of the testing battery rather than measurement error alone.

To further evaluate the contribution of each test to the overall battery structure, the item-rest correlations were analysed. The highest coefficient was observed for Test A (0.302), indicating an acceptable relationship with the remaining tests and suggesting that this component is the most representative of the common mathematical competence assessed by the battery.

In contrast, Test B (0.086) and Test C (0.078) demonstrated very weak relationships with the rest of the battery. This may indicate that these tests capture more specific procedural or symbolic micro-skills, such as isolated algebraic transformations or coefficient handling, rather than the broader determinant-solving competence measured by the other tasks. A likely reason for this pattern is the pronounced ceiling effect, which substantially reduces score variability and weakens correlation coefficients.

Test D (0.168) showed a low but interpretable item-rest correlation, suggesting that it functions as a partially integrative task while still reflecting a somewhat distinct stage of reasoning within the stepwise solution process.

The item difficulty analysis further confirmed the presence of a strong ceiling effect. Since all test scores were normalized to the [0,1] interval, the difficulty index numerically coincided with the mean score for each test. The obtained values were 0.974 for Test A, 0.948 for Test B, 0.930 for Test C, and 0.857 for Test D, while the normalized overall score reached 0.904. According to Classical Test Theory, these values indicate that all four tests can be classified as relatively easy tasks, with Test D being relatively the most difficult component of the battery. Such high difficulty indices explain the reduced score variance, pronounced ceiling effect, and relatively low corrected item-total correlations observed in the previous analyses.

Overall, the results show that the proposed stepwise tasks are completed successfully by most students, although the four-task structure remains psychometrically heterogeneous at the overall level. To better understand how stepwise decomposition affects psychometric coherence at the task level, a more detailed analysis of Task C was conducted.

*4.3 Task C as a Representative Example of Stepwise Task Analysis*

The results of the detailed psychometric analysis of Task C are presented in Figure 10 and Tables 3,4.

Task C consists of five sequential substeps corresponding to the determinant expansion and the subsequent solution of the quadratic equation: identification of the coefficients $x^2$, $x$, and the constant term, followed by determination of the roots $x_1$ and $x_2$.

**Table 3.** ***Descriptive Statistics for Test C***

| | **x^2** | **x** | **1** | **x1** | **x2** | **Score** |
|---|---|---|---|---|---|---|
| **Mean** | 0.959 | 0.917 | 0.925 | 0.934 | 0.917 | 0.930 |
| **Standard deviation** | 0.167 | 0.256 | 0.242 | 0.226 | 0.256 | 0.172 |
| **Skewness** | -5.07 | -3.19 | -3.44 | -3.72 | -3.19 | -3.68 |
| **Kurtosis** | 26.6 | 8.81 | 10.6 | 12.8 | 8.81 | 15.4 |
| **Item-rest correlation** | 0.501 | 0.550 | 0.471 | 0.747 | 0.645 | - |

All five substeps of Task C yielded high mean scores, ranging from 0.917 to 0.959, while the normalized total score for Task C reached 0.930. Since all substep scores were normalized to the [0;1] interval, these mean values are numerically equivalent to the difficulty indices, indicating that all stages of the task can be classified as easy for the examined sample. The highest success rate was observed for the identification of the $x^2$ coefficient ($p = 0.959$), whereas slightly lower values were obtained for the $x$ coefficient and the $x_2$ root ($p = 0.917$), suggesting that these stages require more careful algebraic manipulation.

The comparison of the mean, median, and mode values confirms the presence of a pronounced ceiling effect, as all central tendency indicators reached or approached the maximum possible score. This interpretation is further supported by the strongly negative skewness coefficients (from −5.07 to −3.19) and high kurtosis values (from 8.81 to 26.6), which indicate a highly concentrated score distribution near perfect performance.

**Table 4.** ***Scale Reliability Statistics***

| | Mean | SD | Cronbach's α | McDonald's ω |
|---|---|---|---|---|
| **scale** | 0.930 | 0.172 | 0.794 | 0.808 |

The reliability analysis demonstrated good internal consistency of Task C. The obtained Cronbach's alpha was 0.794, while McDonald's omega reached 0.808, confirming that the five sequential substeps form a sufficiently homogeneous psychometric structure. This result indicates that the decomposition of determinant-based equation solving into micro-operations successfully preserves the internal coherence of the mathematical procedure.

The corrected item-rest correlation coefficients ranged from 0.471 to 0.747, indicating that all substeps contribute meaningfully to the total structure of the task. The strongest contribution was observed for the determination of the first root $x_1$ (0.747), followed by the second root $x_2$ (0.645), which suggests that the root-finding stage is the most representative part of the entire determinant-solving process.

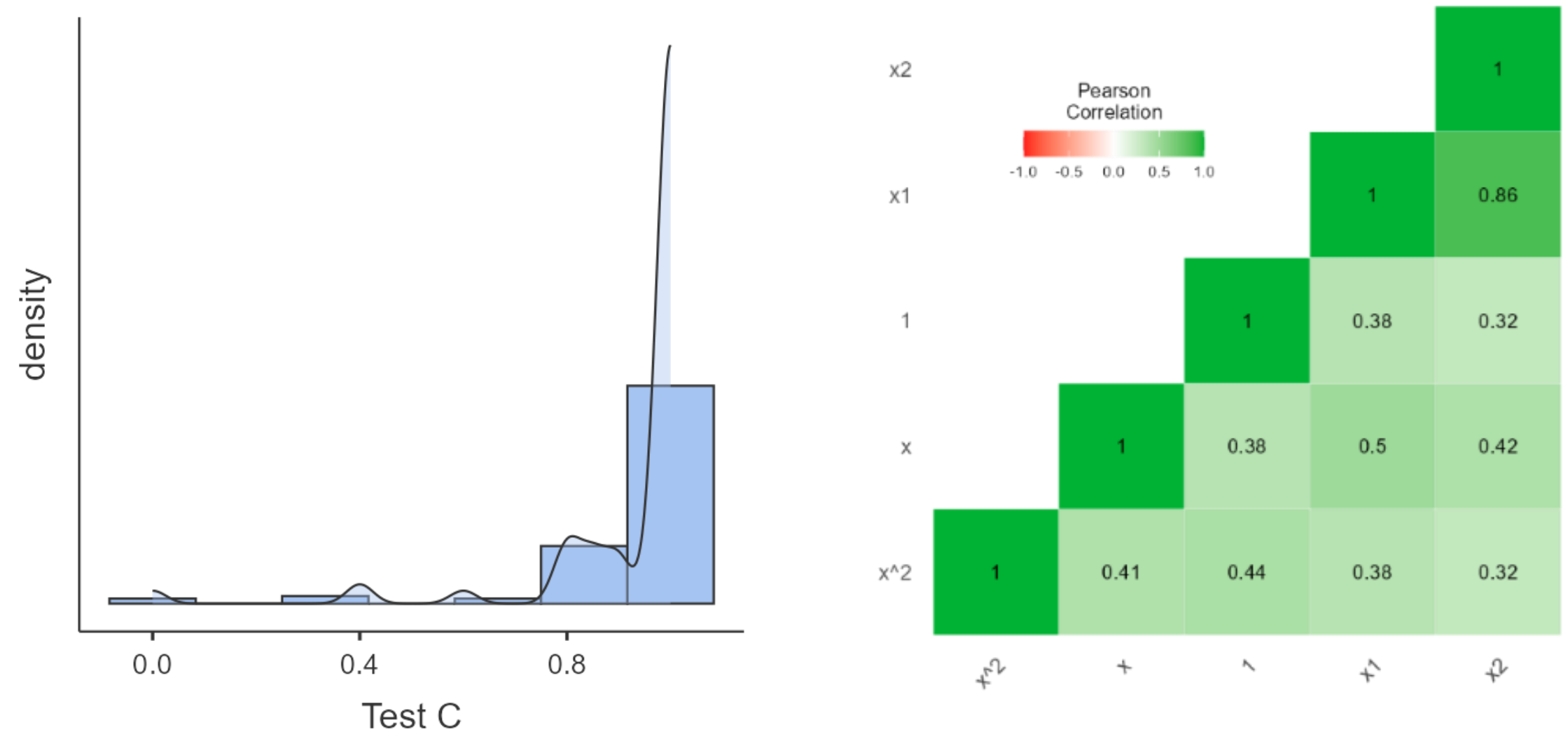

(*a*) Histogram (b) Correlation heatmap

**Figure 10:** Histogram and correlation heatmap between substeps of Task C.

The correlations presented in Figure 10.(*b*) further support the logical sequential nature of the solution process. In particular, the strong correlation between $x_1$ and $x_2$ ($r = 0.86$) confirms the close mathematical dependence of both roots within the same quadratic equation. Moderate correlations between the coefficient-identification stages and the root-determination stages ($r = 0.32 - 0.50$) indicate that successful completion of the task requires both correct determinant expansion and subsequent algebraic transformations.

Overall, these findings show that Task C provides a psychometrically coherent example of an effective stepwise assessment format. The combination of high reliability, acceptable corrected item-total correlations, and interpretable difficulty indices demonstrates that the proposed decomposition of determinant-based equation solving into sequential Moodle-based micro-steps provides a reliable way to assess students' mathematical reasoning.

## 5. Discussions

The contrast between the overall four-task battery and the detailed behaviour of Task C reveals an important feature of stepwise digital assessment in higher mathematics.

At the level of the four-test battery, the results revealed limited internal homogeneity, as reflected in the relatively low values of Cronbach's alpha and McDonald's omega. This suggests that the four tests capture related but distinct components of mathematical competence, including determinant expansion, algebraic transformations, and sequential problem-solving operations. From the perspective of mathematics education, this heterogeneity is expected rather than problematic, since complex competence naturally consists of partially independent micro-skills.

At the same time, the detailed analysis of Task C, chosen as an illustrative example of a typical compact stepwise mathematical procedure, demonstrated that when a mathematical procedure is decomposed into logically dependent micro-operations, the resulting task structure becomes psychometrically coherent. The acceptable reliability coefficients, strong corrected item-rest correlations, and interpretable sequential correlation structure indicate that the internal steps of Task C successfully preserve the logical flow of determinant-based equation solving.

Importantly, the present study should be viewed not as a report on a single classroom implementation, but as a methodological exemplar demonstrating the full cycle of stepwise mathematical task design, Moodle-based deployment, and psychometric validation. The determinant-based task was intentionally selected as a representative medium-complexity example that makes each stage of this workflow transparent and reproducible.

This observation supports the central idea of the study that digital assessment is most effective when task decomposition follows the natural logic of mathematical reasoning. In other words, the effectiveness of stepwise assessment depends not only on digital implementation but also on the extent to which the decomposition corresponds to authentic stages of mathematical thinking.

This interpretation also aligns with recent research in mathematics e-assessment, where high-quality task design, follow-up error handling, and meaningful intermediate feedback are identified as key conditions for effective digital mathematical problem solving (Gallaun et al., 2022; Kinnear et al., 2024). In this sense, the proposed Moodle-based methodology complements specialised systems such as STACK and Numbas by demonstrating that psychometrically coherent stepwise reasoning tasks can also be effectively implemented within a widely adopted LMS environment.

Among the analysed stepwise tasks, Task C and structurally similar multistage tasks appear to be particularly suitable candidates for further STACK-based implementation in Moodle. Unlike conventional numerical-response tasks, the sequential structure of Task C naturally aligns with STACK's capabilities for symbolic input validation, algebraic equivalence checking, and stepwise feedback generation. The decomposition into coefficient identification and root determination can be directly translated into STACK answer trees, enabling automatic verification not only of final answers but also of intermediate symbolic transformations.

A limitation of the present study is that the detailed psychometric analysis was demonstrated on one representative determinant-based task and a single student cohort. Although this task was intentionally selected as a methodological exemplar of a compact multistage mathematical procedure, the generalizability of the reported psychometric patterns to other mathematical domains and larger student populations requires further empirical validation. Future studies should therefore extend this analysis to broader families of stepwise tasks, multiple mathematical topics, and fully STACK-integrated symbolic task variants.

The transition of such tasks to STACK may further strengthen the proposed methodology through more precise symbolic validation, adaptive stepwise feedback, and parameterized task generation.

More broadly, the present results suggest that the future development of digital mathematics assessment should focus on task-level psychometric optimization rather than relying exclusively on aggregate indicators of heterogeneous task batteries. The contrast between the four-test battery and Task C demonstrates that mathematically coherent micro-step design may provide a stronger foundation for reliable digital assessment than broader heterogeneous task collections.

Taken together, these findings suggest that psychometrically optimized micro-step task design may provide a strong methodological basis for the broader integration of symbolic STACK-supported assessment in Moodle.

## Conclusion

This study presents a structured approach to designing step-by-step mathematics assessment tasks in a standard Moodle environment. By combining checking of intermediate steps, interactive task formats, and automated feedback, the proposed methodology goes beyond the assessment of final answers only.

The results show that such tasks can be analysed effectively using Classical Test Theory. The reliability indicators suggest that these assessment formats can work as coherent and psychometrically stable tools. In addition, they give a more detailed picture of students' solution processes, which may help provide better feedback and support error diagnosis.

From a practical point of view, the study shows that widely used learning management systems can support more process-oriented mathematics assessment without the need for specialised external platforms. This makes the approach relevant for distance learning and other disrupted educational contexts.

At the same time, the study mainly focuses on task design and psychometric analysis and does not include a direct comparison of learning outcomes across different assessment formats. Future research should therefore examine how step-by-step and interactive tasks affect conceptual understanding, problem-solving performance, and long-term learning.

Overall, the findings suggest that step-by-step digital assessment is a promising way to develop more transparent, structured, and pedagogically meaningful mathematics assessment in digital learning environments. The findings also suggest that this methodology may be extended to STACK-based symbolic and adaptive assessment formats in Moodle, with the potential to improve both the scalability of task generation and the quality of automated feedback.

**Data Availability Statement.**

The dataset underlying the analyses reported in Section 4 is openly available in Zenodo (DOI: 10.5281/zenodo.19480465).

**The conflict-of-interest statement.** This work does not have any conflicts of interest.